\documentclass[12pt,a4paper,reqno]{amsart}
\usepackage{amsmath}
\usepackage{amsfonts}
\usepackage{amssymb}
\usepackage{latexsym}
\usepackage{graphicx}
\usepackage{color}
\usepackage{amscd}
\usepackage{epsfig}
\usepackage{cite}
\newtheorem{problem}{Problem}[section]
\newtheorem{definition}[problem]{Definition}
\newtheorem{lemma}[problem]{Lemma}
\newtheorem{theorem}[problem]{Theorem}

\newtheorem{corollary}[problem]{Corollary}

\title{Generic exponential sums associated to Laurent polynomials in one variable}
\author{Chunlei Liu and Niu Chuanze}
\address{The School of Mathematical Sciences, Beijing Normal University, Beijing 100875}
\thanks{This project is supported by NSFC Grant No.
10671015.}
\begin{document}
\maketitle
{\bf Abstract. }The generic Newton polygons for $L$-functions of
exponential sums associated to Laurent polynomials in one variable
are determined when $p$ is large. The corresponding Hasse
polynomials are also determined.



\section{Introduction}
 \hskip 0.2in
We shall determine the generic Newton polygon of $L$-functions of
exponential sums associated to Laurent polynomials in one variable.

Throughout this paper, $p$ denotes a prime number, and $q$ denotes a
power of $p$. Write $\mathbb{F}_p:=\mathbb{Z}/p\mathbb{Z}$. Let
$\bar{\mathbb{F}}_p$ be a fixed algebraic closure of the finite
field $\mathbb{F}_p$, and $\mathbb{F}_q$ the finite field with $q$
elements in $\bar{\mathbb{F}}_p$.

Let $f$ be a Laurent polynomial over $\mathbb{F}_q$. We assume that
the leading exponents of $f$ are prime to $p$. One associates to $f$
the Artin-Schreier curve
$$C_f:y^q-y=f(x).$$
Let $N_k$ be the number of $\mathbb{F}_{q^k}$-rational points
including the infinities on $C_f$. The zeta function of $C_f$ is
defined by
$$Z(t,C_f,\mathbb{F}_q)=\exp(\sum\limits_{k=1}^{+\infty}N_k\frac{t^k}{k}).$$

Let $\overline{\mathbb{Q}}$ be a fixed algebraic closure of
$\mathbb{Q}$. Let $\psi$ denote any nontrivial character of
$\mathbb{F}_p$ into $\overline{\mathbb{Q}}^{\times}$. Let $V_f$ be
the affine line $\mathbb{A}$ over $\mathbb{F}_q$ if $f$ is a
polynomial, and let $V_f$ be the one-dimensional torus $\mathbb{T}$
over $\mathbb{F}_q$ if $f$ is not a polynomial. We have
$$N_k=q^k+1+\sum\limits_{\alpha\in\mathbb{F}_q^{\times}}S(k,\alpha f,\mathbb{F}_q),$$
where the exponential sum $S(k,f,\mathbb{F}_q)$ is defined by
$$S(k,f,\mathbb{F}_q)=\sum\limits_{x\in V_f(\mathbb{F}_{q^k})}
\psi(\text{Tr}_{\mathbb{F}_{q^k}/\mathbb{F}_p}(f(x))).$$ So we have
$$(1-t)(1-qt)Z(t,C_f,\mathbb{F}_q)=\prod\limits_{\alpha\in\mathbb{F}_q^{\times}}L(t,\alpha f,\mathbb{F}_q),$$
where the $L$-function $L(t,f,\mathbb{F}_q)$ is defined by
$$L(t,f,\mathbb{F}_q)=\exp(\sum\limits_{k=1}^{+\infty}S(k,f,\mathbb{F}_q)\frac{t^k}{k}).$$
It is well-known that the function $L(t,f,\mathbb{F}_q)$ is a
polynomial in $t$ with coefficients in $\overline{\mathbb{Q}}$.

Let $\mathbb{Z}_p$ be the ring of $p$-adic integers, and
$\mathbb{Q}_p$ the field of $p$-adic numbers. Fix an embedding of
$\overline{\mathbb{Q}}$ into $\overline{\mathbb{Q}}_p$. Let
$\text{ord}_p(\cdot)$ be the $p$-adic order function of
$\overline{\mathbb{Q}}_p$, and define the $q$-adic order function as
$\text{ord}_q(\cdot)=\frac{1}{\text{ord}_p(q)}\text{ord}_p(\cdot)$.
As $L(t,f,\mathbb{F}_q)$ has coefficients in
$\overline{\mathbb{Q}}$, one can talk about the $p$-adic absolute
values of its reciprocal roots. These $p$-adic absolute values are
completely determined by the Newton polygon of $L(t,f,\mathbb{F}_q)$
defined as follows.

\begin{definition}Let
$g(t)=1+\sum\limits_{i=1}^uc_it^i$ be a polynomial in $t$ with
coefficients $c_i\in\overline{\mathbb{Q}}_p$. The $q$-adic Newton
polygon of $g$ is the lower convex closure of the points
$$(0,0),(n,\text{ord}_q(c_n)),\ n=1,\cdots,u.$$
\end{definition}

It is very hard to determine the Newton polygon of
$L(t,f,\mathbb{F}_q)$ in general. However, it is easier to give a
good lower bound. The simplest one is the Hodge polygon defined as
follows.
\begin{definition}Let $d$ be a positive integer.
The Hodge polygon of the interval $[0,d]$ is the polygon whose vertices are
$$(n,\frac{n(n+1)}{2d}),\
n=0,1,2,\cdots,d-1.$$ \end{definition}

\begin{definition}Let $d$ and $e$ be positive integers.
The Hodge polygon of $[-e,d]$ is the polygon with initial point
$(0,0)$, end point $(d+e,\frac{d+e}{2})$, and the vertices
$(m+n+1,\frac{m(m+1)}{2e}+\frac{n(n+1)}{2d})$ with $(m,n)$ running
over pairs satisfying
$$-\frac{1}{e}<\frac{m}{e}-\frac{n}{d}<\frac{1}{d},0\leq m<e,0\leq n<d.$$
\end{definition}

Let $\Delta(f)$ be the smallest closed interval of the real line
containing $0$ and the exponents of the monomials of $f$. So
$\Delta(f)=[0,d]$ if $f$ is a polynomial of degree $d$, and
$\Delta(f)=[-e,d]$ if $f$ is a Laurent polynomial with leading term
$a_{-e}x^{-e}+a_dx^d$.  The well-known Hodge bound is stated as the
following theorem.
\begin{theorem}[Hodge bound]The $q$-adic Newton polygon of $L(t,f,\mathbb{F}_q)$
lies above the Hodge polygon of $\Delta(f)$. Moreover, both polygons
have the same initial point and the same end point.\end{theorem}

By Grothendieck's specialization lemma (Confer \cite{K76} and
\cite{W04}), the $q$-adic Newton polygon of $L(t,f,\mathbb{F}_q)$ is
constant for a generic $f$ with fixed $\Delta(f)=\Delta$. That
constant polygon is called the generic Newton polygon of $\Delta$.
Let $d>0$ and $e\geq0$ be integers. Let $D=d$ if $e=0$, and let $D$
be the least common multiple of $d$ and $e$ if $e>0$. We assume that
$D$ is prime to $p$. The following theorem says that in nice
situations the generic polygon coincides with the Hodge polygon.

\begin{theorem}[Stickelberger's theorem \cite{W93}] The generic Newton polygon
of $[-e,,d]$ coincides with its Hodge polygon if and only if
$p\equiv 1(\mod D)$.\end{theorem}

As a special case of a conjecture of Wan \cite{W04}, the following
theorem says that the Hodge bound is approximately the best.
\begin{theorem}[Zhu \cite{Zh03,Zh04,Zhu04}] The generic Newton polygon of $[-e,d]$ goes to its Hodge
polygon as $p$ goes to infinity.\end{theorem}

In proving the above theorem in the case $e=0$, Zhu used Dwork's
$p$-adic theory, a kind of Di\'{e}donne-Manin diagonalization, and
some force computations to produce a list of polynomials she denoted
as $f_t$'s. She then used a kind of maximal-monomial-locating
technique to prove that one of these $f_t$'s does not vanish.
Blache-F\'{e}rard \cite{BF} discovered that Zhu's
maximal-monomial-locating technique can prove the nonvanishing of
$f_0$. This enabled them to get the following theorem.

\begin{theorem}[Blache-F\'{e}rard]
If $p\geq3D$, the generic Newton polygon of $[0,d]$ is the polygon
with vertices
$$(n,\frac{1}{p-1}\sum\limits_{i=1}^n\lceil\frac{pi-n}{d}\rceil),\ n=0,1,\cdots,d-1.$$
\end{theorem}

The condition $p\geq3D$ in the theorem is very clean. To achieve
that clean condition Blache-F\'{e}rard abolished Zhu's
Di\'{e}donne-Manin diagonalization technique, and made recourse to
Dwork's original method.

It should be mentioned that Yang \cite{Ya} computed the Newton
polygons for $L$-function of exponential sums associated to
polynomials of the form $x^d+\lambda x$, and Hong \cite{H01, H02}
computed the Newton polygons for $L$-function of exponential sums
associated to polynomials of degree $4$ and $6$.

From now on we assume that $e>0$. We shall determine the generic
Newton polygon of $[-e,d]$.

Write
$$p_{[0,d]}(n)=\frac{1}{p-1}\sum\limits_{i=1}^n\lceil\frac{pi-n}{d}\rceil,n=0,1,\cdots,d.$$
And write
$$p_{[-e,d]}(0)=0,p_{[-e,d]}(d+e)=\frac{d+e}{2},$$
$$p_{[-e,d]}(k)=\min\limits_{(m,n)\in I_k}\{p_{[0,e]}(m)+p_{[0,d]}(n)\},\
k=1,\cdots,d+e-1,$$where
$$I_k=\{(m,n)\mid m+n+1=k,
-\frac{1}{e}\leq\frac{m}{e}-\frac{n}{d}\leq\frac{1}{d}, 0\leq
m<e,0\leq n<d\}.$$

\begin{definition}\label{arith-polygon}
The arithmetic polygon of $[-e,d]$ is defined to be the graph of the
function on $[0,d+e]$ which is linear between consecutive integers
and takes on the value $p_{[-e,d]}(k)$ at integers
$k=0,1,\cdots,d+e$.\end{definition}

We shall prove the following theorem.
\begin{theorem}\label{generic-newton}The generic Newton polygon of $[-e,d]$ coincides
with its arithmetic polygon if $p\geq3D$.
\end{theorem}
It would be interesting if one can extend the result to twisted
exponential sums and to exponential sums associated to functions
studied in \cite{Zhu04} and \cite{BFZ}.

{\bf Acknowledgement. }The first author thanks Lei Fu and Daqing Wan
for discussions. He also thanks Daqing Wan for sending the work of
Blache and F\'{e}rard.

\section{The arithmetic polygon}
Recall that, for $k=1,\cdots,d+e-1$,
$$I_k=\{(m,n)\mid m+n+1=k,
-\frac{1}{e}\leq\frac{m}{e}-\frac{n}{d}\leq\frac{1}{d}, 0\leq
m<e,0\leq n<d\}.$$ Let $V_k$ be the subset of $I_k$ consisting pairs
at which the function
$$(m,n)\mapsto p_{[0,e]}(m)+p_{[0,d]}(n)$$
takes on the minimal value. In this section we shall prove the
following theorem.\begin{theorem}\label{arith-polygon2}Let $p>3D$.
Then the arithmetic polygon of $[-e,d]$ is convex. Moreover,
$(k,p_{[-e,d]}(k))$ ($0<k<d+e$) is a vertex if and only if $V_k$
contains only one pair.
\end{theorem}

We begin with the following lemma.
\begin{lemma}The set
$I_k$ contains one or two pairs. If $I_k=\{(m,n)\}$, then
$$-\frac{1}{e}<\frac{m}{e}-\frac{n}{d}<\frac{1}{d}.$$
 If
$I_k$ contains exactly two pairs, then it is of form
$\{(m,n),(m+1,n-1)\}$ with
$$\frac{m+1}{e}=\frac{n}{d}.$$
\end{lemma}
{\it Proof. }Define a degree function on $\mathbb{Z}$ by
$$\deg(i)=\left\{
            \begin{array}{ll}
              i/d, & \hbox{ }i\geq 0, \\
              -i/e, & \hbox{ }i\leq 0.
            \end{array}
          \right.
$$There is a
positive integer $u$ such that
$$k=\#\{i\in \mathbb{Z}\mid \deg(i)\leq u/D\},$$ or
$$\#\{i\in \mathbb{Z}\mid \deg(i)\leq u/D\}<k<\#\{i\in\mathbb{Z}\mid \deg(i)\leq (u+1)/D\}.$$

If $k=\#\{i\in \mathbb{Z}\mid \deg(i)\leq u/D\}$, then $I_k$ is of
form $\{(m,n)\}$ with
$$-\frac{1}{e}<\frac{m}{e}-\frac{n}{d}<\frac{1}{d}.$$
If $\#\{i\in \mathbb{Z}\mid \deg(i)\leq
u/D\}<k<\#\{i\in\mathbb{Z}\mid \deg(i)\leq (u+1)/D\}$, then $I_k$ is
of form $\{(m,n),(m+1,n-1)\}$ with
$$\frac{m+1}{e}=\frac{n}{d}.$$ The lemma is proved.

It is easy to see that Theorem \ref{arith-polygon2} follows from the
following three theorems.

\begin{theorem}\label{arith-polygon3}Let $k=1,2,\cdots,d+e-1$. If $V_k$ contains two pairs,
then
$$2p_{[-e,d]}(k)=p_{[-e,d]}(k-1)+p_{[-e,d]}(k+1).$$
\end{theorem}

\begin{theorem}\label{arith-polygon4}Let $k=1,2,\cdots,d+e-1$. If $I_k$ contains two pairs
but $V_k$ contains only one pair, then
$$2p_{[-e,d]}(k)<p_{[-e,d]}(k-1)+p_{[-e,d]}(k+1).$$
\end{theorem}

\begin{theorem}\label{arith-polygon5}Let $p>3D$. Let $k=1,2,\cdots,d+e-1$. If $I_k$ contains only one pair,
then
$$2p_{[-e,d]}(k)<p_{[-e,d]}(k-1)+p_{[-e,d]}(k+1).$$
\end{theorem}
{\it Proof of Theorem \ref{arith-polygon3}.} Suppose that $V_k$
contains two pairs. Then so does $I_k$. Assume that
$I_k=\{(m,n),(m+1,n-1)\}$. Then $\frac{m+1}{e}=\frac{n}{d}$. It
follows that $I_{k-1}=\{(m,n-1)\}$ and $I_{k+1}=\{(m+1,n)\}$. Note
that
$$p_{[-e,d]}(k)=p_{[0,e]}(m)+p_{[0,d]}(n)=p_{[0,e]}(m+1)+p_{[0,d]}(n-1).$$
It follows
that
$$2p_{[-e,d]}(k)=p_{[0,e]}(m)+p_{[0,d]}(n)+p_{[0,e]}(m+1)+p_{[0,d]}(n-1)
=p_{[-e,d]}(k-1)+p_{[-e,d]}(k+1).$$
Theorem \ref{arith-polygon3} is proved.

{\it Proof of Theorem \ref{arith-polygon4}.} Assume that
$I_k=\{(m,n),(m+1,n-1)\}$. Then $\frac{m+1}{e}=\frac{n}{d}$. It
follows that $I_{k-1}=\{(m,n-1)\}$ and $I_{k+1}=\{(m+1,n)\}$.
Without loss of generality, we assume that $V_k=\{(m,n)\}$. Then
$$p_{[0,e]}(m)+p_{[0,d]}(n)<p_{[0,e]}(m+1)+p_{[0,d]}(n-1).$$
It follows
that$$2p_{[0,e]}(m)+2p_{[0,d]}(n)<p_{[0,e]}(m)+p_{[0,d]}(n)+p_{[0,e]}(m+1)+p_{[0,d]}(n-1).$$
That is,
$$2p_{[-e,d]}(k)<p_{[-e,d]}(k-1)+p_{[-e,d]}(k+1).$$
Theorem \ref{arith-polygon4} is proved.

{\it Proof of Theorem \ref{arith-polygon5}.} Assume that
$I_k=\{(m,n)\}$. Then
$$-\frac{1}{e}<\frac{m}{e}-\frac{n}{d}<\frac{1}{d}.$$
Let $(m_1,n_1)\in V_{k-1}$. Then $m_1=m$ or $n_1=n$. Without loss of
generality, we assume that $m_1=m$. Then $n_1=n-1$. Let
$(m_2,n_2)\in V_{k+1}$. Then $m_2=m$ or $n_2=n$.

First, we assume that $m_2=m$. Then $n_2=n+1$. Note that
$$p_{[0,d]}(n+1)-p_{[0,d]}(n)\geq\frac{1}{p-1}(\lceil(p-1)\frac{n+1}{d}\rceil-1),$$
$$p_{[0,d]}(n)-p_{[0,d]}(n-1)\leq\frac{1}{p-1}\lceil(p-1)\frac{n}{d}\rceil,$$
and
$$\lceil(p-1)\frac{n}{d}\rceil<\lceil(p-1)\frac{n+1}{d}\rceil-1.$$
It follows that
$$2p_{[0,d]}(n)<p_{[0,d]}(n+1)+p_{[0,d]}(n-1).$$
Therefore
$$2p_{[0,e]}(m)+2p_{[0,d]}(n)=p_{[0,e]}(m)+p_{[0,d]}(n+1)+p_{[0,e]}(m)+p_{[0,d]}(n-1).$$
That is,
$$2p_{[-e,d]}(k)<p_{[-e,d]}(k-1)+p_{[-e,d]}(k+1).$$

Secondly, we assume that $n_2=n$. Then $m_2=m+1$. Note that
$$p_{[0,e]}(m+1)-p_{[0,e]}(m)\geq\frac{1}{p-1}(\lceil(p-1)\frac{m+1}{e}\rceil-1),$$
$$p_{[0,d]}(n)-p_{[0,d]}(n-1)\leq\frac{1}{p-1}\lceil(p-1)\frac{n}{d}\rceil,$$
and
$$\lceil(p-1)\frac{n}{d}\rceil<\lceil(p-1)\frac{m+1}{e}\rceil-1.$$
It follows that
$$p_{[0,e]}(m)+p_{[0,d]}(n)<p_{[0,e]}(m+1)+p_{[0,d]}(n-1).$$
Therefore
$$2p_{[0,e]}(m)+2p_{[0,d]}(n)<p_{[0,e]}(m+1)+p_{[0,d]}(n)+p_{[0,e]}(m)+p_{[0,d]}(n-1).$$
That is,
$$2p_{[-e,d]}(k)<p_{[-e,d]}(k-1)+p_{[-e,d]}(k+1).$$
Theorem \ref{arith-polygon5} is proved.

\section{Hasse polynomial}
For $\vec{a}=(a_{-e},\cdots,a_d)\in\mathbb{F}_q^{d+e+1}$, we write
$$f_{\vec{a}}(x)=\sum\limits_{i=-e}^da_ix^i.$$
It is easy to see that the Newton polygon of
$L(t,f_{\vec{a}},\mathbb{F}_q)$ is independent of $a_0$. So one can
take $a_0$ to be any preferred number. We take $a_0=1$ so that
Lemmas \ref{artin-hasse-plus-main-term} and
\ref{artin-hasse-minus-main-term} are expressed in a simpler form.

In this section we define a polynomial $H$ such that the Newton
polygon of $L(t,f_{\vec{a}},\mathbb{F}_q)$ coincides with the
generic Newton polygon of $[-e,d]$ if and only if $H(\vec{a})\neq0$.

\begin{definition}Let $k=1,2,\cdots,d+e-1$ be such that $V_k=\{(m,n)\}$.
 We define $S_k$ to be the set of
permutations $\tau$ of $\{-m,-m+1,\cdots,n\}$ such that
$$\tau(i)\left\{
                                                        \begin{array}{ll}
\geq n-d\{-\frac{pi-n}{d}\}, & \hbox{ if } i>0,\\
=0, & \hbox{ if } i=0,\\
\leq -m+e\{\frac{pi+m}{e}\}, & \hbox{ if }i<0.
                                                        \end{array}
                                                      \right.
$$
\end{definition}

Let
$$E(t)=\exp(\sum_{i=0}^{\infty}\frac{t^{p^i}}{p^i}).$$
It is a power series in $\mathbb{Z}_p[[t]]$, and we call it the
Artin-Hasse exponential series. We write
$$E(t)=\sum\limits_{n=0}^{+\infty}\lambda_nt^n.$$
\begin{definition}Let $k=1,2,\cdots,d+e-1$ be such that $V_k=\{(m,n)\}$. We write
$$r_i=\left\{
                                                   \begin{array}{ll}
                                                            n-d\{-\frac{pi-n}{d}\}+d, & \hbox{ } 1\leq i\leq n,\\
                                                        m-e\{\frac{pi+m}{e}\}+e, & \hbox{
}-m\leq i\leq -1.
                                                          \end{array}
                                                        \right.
$$
We define a polynomial $H_k$ in the variables $x_{-e},\cdots,x_d,$
by
$$H_k(\vec{x})=\sum\limits_{\tau\in S_k}u_{\tau} \prod_{i=-m}^{-1}x_{-r_i-\tau(i)}
\prod_{i=1}^nx_{r_i-\tau(i)},$$ where
$$
u_{\tau}=\text{sgn}(\tau)(\prod_{i=1}^n\lambda_{\lfloor\frac{pi-\tau(i)}{d}\rfloor}
\lambda_{\lceil\{\frac{pi-\tau(i)}{d}\}\rceil})\prod_{i=-m}^{-1}\lambda_{\lfloor\frac{-pi+\tau(i)}{e}\rfloor}
\lambda_{\lceil\{\frac{-pi+\tau(i)}{e}\}\rceil}\in\mathbb{Z}_p^{\times}.$$
\end{definition}

\begin{definition}The Hasse polynomial $H$ of $[-e,d]$
is defined by $$H=x_dx_{-e}\prod_{\#V_k=1}\bar{H}_k,$$ where
$\bar{H}_k$ is the reduction of $H_k$ modulo $p$.
\end{definition}

We shall prove the following theorem.
\begin{theorem}\label{non-zero}
The Hasse polynomial $H$ of $[-e,d]$ is non-zero.
\end{theorem}

Let $k=1,2,\cdots,d+e-1$ be such that $V_k=\{(m,n)\}$. It is easy to
see that Theorem \ref{non-zero} follows from the following one.

\begin{theorem}\label{multiplicity-one}Among
the monomials
$$\prod_{i=-m}^{-1}x_{-r_i-\tau(i)}
\prod_{i=1}^nx_{r_i-\tau(i)},\ \tau\in S_k,
$$ there is a monomial which appears exactly once.
\end{theorem}

That theorem plays a crucial role in the determination of the
generic Newton polgon of $[-e,d]$. In the case $e=0$,
Blache-F\'{e}rard \cite{BF} used Zhu's maximal-monomial-locating
technique to prove the theorem. In the case $e>0$, the
maximal-monomial-locating technique no longer works. Fortunately, a
minimal-monomial-locating technique will play the role.

Set $x_1<x_2<\cdots<x_d$ and $x_{-1}<x_{-2}<\cdots<x_{-e}$. Define
$\prod_{i\in I}x_i>\prod_{j\in J}x_j$ and $\prod_{i\in
I}x_{-i}>\prod_{j\in J}x_{-j}$ if $I$ and $J$ are finite subsets of
positive integers and there is an $i\in I$ which is greater than all
$j\in J$. Define $g_1g_3\geq g_2g_4$ if $g_1,g_2,g_3,g_4$ are
monomials such that $g_1\geq g_2$ and $g_3\geq g_4$.

It is easy to see that Theorem \ref{multiplicity-one} follows from
the following one.
\begin{theorem}\label{minimal-monomial}Among
the monomials
$$\prod_{i=-m}^{-1}x_{-r_i-\tau(i)}
\prod_{i=1}^nx_{r_i-\tau(i)},\ \tau\in S_k,
$$ the minimal monomial appears exactly once.
\end{theorem}
{\it Proof. } Note that $r_i\neq r_j$ and $r_{-i}\neq r_{-j}$ if $i$
and $j$ are distinct positive integers. So we can order them such
that
$$r_{i_1}>r_{i_2}>\cdots>r_{i_n},\ i_j>0,$$ and
$$r_{t_1}>r_{t_2}>\cdots>r_{t_m},\ t_j<0.$$
Note that $r_{i_1}\leq n+d$ and $r_{t_1}\leq m+e$. So we have
$$r_{i_j}\leq n+d+1-j,\text{ and } r_{t_j}\leq m+e+1-j.$$
Recall that $\tau\in S_k$ if and only if $\tau(i)\geq r_i-d$ if
$i>0$, and $\tau(i)\leq -r_i+e$ if $i<0$. Hence, if we define
$\tau_0$ by
$$\tau_0(i_j)=n+1-j, \text{ and }\tau_0(t_j)=-(m+1-j),$$
then $\tau_0\in S_k$.

We claim that, for any $\tau\in S_k$,
$$\prod_{j=1}^nx_{r_{i_j}-\tau(i_j)}\geq\prod_{j=1}^nx_{r_{i_j}-(n+1-j)}
$$
with equality holding if and only if $\tau(i_j)=n+1-j$ for all
$1\leq j\leq n$. Suppose that $\tau(i_j)\neq n+1-j$ for some $1\leq
j\leq n$. Let $j_0$ be the least one with this property. Then
$\tau(i_{j_0})<n+1-j_0$. Hence
$$r_{i_{j_0}}-\tau(i_{j_0})>r_{i_{j_0}}-(n+1-j_0)\geq r_{i_j}-(n+1-j),\text{ for all }j\geq j_0 .$$
Therefore
$$\prod_{j=1}^nx_{r_{i_j}-\tau(i_j)}>\prod_{j=1}^nx_{r_{i_j}-(n+1-j)}
$$
as claimed.

Similarly, we can prove that, for any $\tau\in S_k$,
$$\prod_{j=1}^mx_{-r_{t_j}-\tau(t_j)}\geq\prod_{j=1}^mx_{-r_{t_j}+(m+1-j)}$$
with equality holding if and only if $\tau(t_j)=-(m+1-j)$ for all
$1\leq j\leq m$. It follows that the monomial
$$\prod_{j=1}^nx_{r_{i_j}-(n+1-j)}\prod_{j=1}^mx_{-r_{t_j}+(m+1-j)}$$
is minimal and occurs in the monomials
$$\prod_{i=-m}^{-1}x_{-r_i-\tau(i)}
\prod_{i=1}^nx_{r_i-\tau(i)},\ \tau\in S_k.
$$ The theorem is proved.

\section{Dwork's $p$-adic analytic method}
In this section we give a brief survey on Dwork's $p$-adic analytic
method. Proofs of theorems in this section may be omitted.
Interested readers may consult \cite{Dw62, Dw64} and \cite{AS87,
AS89} for detailed proofs.

 Write $\mathbb{Z}_q:=\mathbb{Z}_p[\mu_{q-1}]$ and
$\mathbb{Q}_q:=\mathbb{Q}_p(\mu_{q-1})$, where $\mu_n$ is the group
of $n$-th roots of unity.

Recall that
$$E(t)=\exp(\sum_{i=0}^{\infty}\frac{t^{p^i}}{p^i})=\sum\limits_{n=0}^{+\infty}\lambda_nt^n \in
\mathbb{Z}_p[[t]]$$ is the Artin-Hasse exponential series. Choose
$\pi\in\mathbb{Q}_p(\mu_p)$ such that $E(\pi)=\psi(1)$. We have
$\text{ord}_p(\pi)=\frac{1}{p-1}$ and
$\sum\limits_{i=0}^{\infty}\frac{\pi^{p^i}}{p^i}=0$.

Let $L$ be the Banach space over $\mathbb{Q}_{q}[\pi^{1/D}]$ with
formal basis $\pi^{\deg(i)}x^i$, $i\in \mathbb{Z}$. That is,
$L=L_0\otimes_{\mathbb{Z}_q}\mathbb{Q}_q$ with
$$L_0=\{\sum_{i\in
\mathbb{Z}}c_i\pi^{\deg(i)}x^i:\
 c_i\in\mathbb{Z}_q[\pi^{1/D}] \}.$$
The space is closed under multiplication. So it is an algebra.

For $\vec{a}=(a_{-e},\cdots,a_d)\in\mathbb{F}_q^{d+e+1}$, we write
$$E_{\vec{a}}(x) :=\prod\limits_{i=-e}^dE(\pi
\hat{a}_ix^i),$$ where $\hat{a}_i$ is the Teichm\"{u}ller lifting of
$a_i$. As each $E(\pi \hat{a}_ix^i)$ lies in $L$, so does
$E_{\vec{a}}$.

The Galois group $\text{Gal}(\mathbb{Q}_q/\mathbb{Q}_p)$ acts on $L$
but keeps $\pi^{1/D}$ and $x$ fixed. Let $\sigma$ be the Frobenius
element of that Galois group. Write
$$\hat{E}_{\vec{a}}(x)=\prod\limits_{j=0}^{+\infty}E_{\vec{a}}^{\sigma^j}(x^{p^j}).$$
Define an operator $\partial:L\rightarrow L$ by
$$\partial(g)=xg'(x)+g(x)x\frac{d\log\hat{E}_{\vec{a}}(x)}{dx}.$$
It is easy to see that $L_0$ is stable under $\partial$.

Define an operator $\psi_p:L\rightarrow L$ by
$$\psi_p(\sum\limits_{i\in \mathbb{Z}} c_ix^i)=\sum\limits_{i\in \mathbb{Z}}
c_{pi}x^i.$$ And write
$$\Psi_p:=\sigma^{-1}\circ\psi_p\circ
E_{\vec{a}}.$$ That is,
$$\Psi_p(g)=\sigma^{-1}(\psi_p(gE_{\vec{a}})).$$
Note that $\Psi_p$ is $\mathbb{Q}_p[\pi^{\frac{1}{D}}]$-linear, but
$\mathbb{Q}_q[\pi^{\frac{1}{D}}]$-semi-linear.

Define $\Psi_{p^n}=\Psi_p^n$. So $\Psi_{q^n}=\Psi_q^n$. It is easy
to check that $\Psi_q$ is $\mathbb{Z}_q[\pi^{1/D}]$-linear.
Moreover, we have
$$q\partial\Psi_q=\Psi_q\partial.$$

Let $\bar{\Psi}_p$ be the induced operator of $\Psi_p$ on
$L/(\partial L)$. We have the following three
theorems.
\begin{theorem}\label{trace-formula} We have
$$L(s,f_{\vec{a}},\mathbb{F}_{q})
=\det(1-s\bar{\Psi}_q\mid L/(\partial L)\text{ over
}\mathbb{Q}_q(\pi^{1/D})).$$
\end{theorem}

\begin{theorem}\label{q2p}The $q$-adic Newton polygons of
 $\det(1-s^b\bar{\Psi}_q\mid L/(\partial L)\text{ over }\mathbb{Q}_q(\pi^{1/D}))$ and
$\det(1-s\bar{\Psi}_p\mid L/(\partial L)\text{ over
}\mathbb{Q}_p(\pi^{1/D}))$ coincide.\end{theorem}

\begin{theorem}Over $\mathbb{Z}_q[\pi^{1/D}]$, the lattice $L_0/(\partial L_0)$ has a
basis represented by
$$\pi^{\deg(i)}x^i,\ -e\leq i\leq d-1.$$\end{theorem}
\section{Elementary estimates}
In this section we give some elementary estimates on the matrix
coefficients of the operator $\bar{\Psi}_p$ on $L/(\partial L)$.

Write
$$E_{\vec{a}}(x)=
\sum\limits_{i\in \mathbb{Z}}\gamma_ix^i.$$ We have
$$\gamma_i=\sum\limits_{\sum\limits_{j=-e}^djn_j=i}\pi^{\sum\limits_{j=-e}^dn_j}
\prod_{j=-e}^d\lambda_{n_j}\hat{a}_j^{n_j}.$$

\begin{definition}We write $\alpha=O(\pi^t)$ to mean that $\text{ord}_{\pi}(\alpha)\geq t$, where
$\text{ord}_{\pi}(\cdot)=\frac{1}{\text{ord}_p(\pi)}\text{ord}_p(\cdot)$.\end{definition}

\begin{lemma}\label{artin-hasse-plus-main-term}If $i\geq0$,
$$\gamma_i=\pi^{\lceil\frac{i}{d}\rceil}\lambda_{\lfloor
\frac{i}{d}\rfloor}\lambda_{\lceil\{\frac{i}{d}\}\rceil}\hat{a}_d^{\lfloor\frac{i}{d}\rfloor}\hat{a}_{d\{\frac{i}{d}\}}
+O(\pi^{\lceil\frac{i}{d}\rceil+1}).$$
\end{lemma}
{\it Proof. }If $\sum\limits_{j=-e}^djn_j=i$ ($n_j\geq0$), then
$\sum\limits_{j=-e}^dn_j\geq\lceil\frac{i}{d}\rceil$ with equality
holding if and only if $$n_j=\left\{
                               \begin{array}{ll}
                                 \lfloor\frac{i}{d}\rfloor, & \hbox{ } j=d\\
                                 \lceil\{\frac{i}{d}\}\rceil, & \hbox{ }j=d\{\frac{i}{d}\} \\
                                 0, & \hbox{ otherwise.}
                               \end{array}
                             \right.
$$
The lemma now follows.

Similarly, we have the following lemma.
\begin{lemma}\label{artin-hasse-minus-main-term}If $i<0$,
$$\gamma_i=\pi^{\lceil\frac{-i}{e}\rceil}\lambda_{\lfloor
\frac{-i}{e}\rfloor}\lambda_{\lceil\{\frac{-i}{e}\}\rceil}\hat{a}_{-e}^{\lfloor\frac{-i}{e}\rfloor}
\hat{a}_{-e\{\frac{-i}{e}\}} +O(\pi^{\lceil\frac{-i}{e}\rceil+1}).$$
\end{lemma}

From the last two lemmas we infer the following corollary.
\begin{corollary}\label{artin-hasse-estimate}We have
$$\gamma_i=O(\pi^{\lceil\deg(i)\rceil}).$$
\end{corollary}

Let $F=(F_{ij})_{-e\leq i,j\leq d-1}$ be the matrix defined by
$$\psi_p\circ E_{\vec{a}}(x^j)\equiv\sum\limits_{i=-e}^{d-1}F_{ij}x^i (\mod \partial L).$$

\begin{lemma}\label{frobenius-main-term}Let $p\geq3D$, and $-e\leq i,j\leq d-1$. We have
$$F_{ij}=\gamma_{pi-j}+\left\{
           \begin{array}{ll}
              O(\pi^{\lfloor\deg(pi)\rfloor+2}), & \hbox{ }i\neq-e\\
            O(\pi^p), & \hbox{ }i=-e.
           \end{array}
         \right.
$$ \end{lemma}
{\it Proof. }We have
$$\psi_p\circ
E_{\vec{a}}(x^j)=\sum\limits_{i_0\in
\mathbb{Z}}\gamma_{pi_0-j}x^{i_0}=\sum\limits_{i_0=-e}^{d-1}\gamma_{pi_0-j}x^{i_0}
+\sum\limits_{i_0\not\in\{-e,\cdots,d-1\}}\gamma_{pi_0-j}x^{i_0}.$$
For $i_0\not\in\{-e,\cdots,d-1\}$, we write
$$\pi^{\deg(i_0)}x^{i_0}=\sum\limits_{i=-e}^{d-1}c_{ii_0}\pi^{\deg(i)}x^i(\mod \partial L),\ c_{ii_0}\in\mathbb{Z}_q[\pi^{1/(de)}].$$
Then$$\psi_p\circ
E_{\vec{a}}(x^j)=\sum\limits_{i=-e}^{d-1}x^i(\gamma_{pi-j}
+\sum\limits_{i_0\not\in\{-e,\cdots,d-1\}}c_{ii_0}\pi^{\deg(i)-\deg(i_0)}\gamma_{pi_0-j})(\mod
\partial L).$$ It follows that
$$F_{ij}=\gamma_{pi-j}
+\sum\limits_{i_0\not\in\{-e,\cdots,d-1\}}c_{ii_0}\pi^{\deg(i)-\deg(i_0)}\gamma_{pi_0-j}.$$
If $i_0\not\in\{-e,-(e-1),\cdots,d-1\}$, and $i\neq-e$, we
have$$\deg(i)-\deg(i_0)+\text{ord}_{\pi}(\gamma_{pi_0-j})
\geq\deg(i)-\deg(i_0)+\deg(pi_0)-1$$$$\geq\deg(pi)+(p-1)(\deg(i_0)-\deg(i))-1
$$$$\geq\lfloor\deg(pi)\rfloor+\frac{p-1}{D}-1\geq\lfloor\deg(pi)\rfloor+2.$$
If $i_0\not\in\{-e,-(e-1),\cdots,d-1,d\}$, and $i=-e$, we also
have$$\deg(i)-\deg(i_0)+\text{ord}_{\pi}(\gamma_{pi_0-j})\geq\lfloor\deg(pi)\rfloor+2.$$
If  $i_0=d$, and $i=-e$,we
have$$\deg(i)-\deg(i_0)+\text{ord}_{\pi}(\gamma_{pi_0-j})\geq p.$$
Therefore$$F_{ij}=\gamma_{pi-j}+\left\{
           \begin{array}{ll}
             O(\pi^{\lfloor\deg(pi)\rfloor+2}), & \hbox{ } i\neq-e\\
             O(\pi^p), & \hbox{ }i=-e.
           \end{array}
         \right.
$$The lemma is proved.
\section{Generic polygon}
In this section we prove Theorem \ref{generic-newton}. It follows
immediately from the following theorem.
\begin{theorem}\label{p-adic-riemann}Let
$p\geq3D$. Then the $q$-adic Newton polygon of
$L(t,f_{\vec{a}},\mathbb{F}_q)$ coincides with the arithmetic
polygon of $[-e,d]$ if and only if $H(\vec{a})\ne0$.
\end{theorem}

Write
$$\det(1-s\bar{\Psi}_p\mid L/(\partial L)\text{ over
}\mathbb{Q}_p(\pi^{1/D}))=\sum\limits_{i=0}^{b(d+e)}(-1)^ic_is^i.$$
By Theorems \ref{q2p} and \ref{arith-polygon2}, Theorem
\ref{p-adic-riemann} follows from the following two theorems.

\begin{theorem}\label{non-vertex}Let $p>3D$. Let $k=1,2,\cdots,d+e-1$ be such that $V_k$ contains two pairs.
Then
$$\text{ord}_q(c_{bk})\geq p_{[-e,d]}(k).$$
\end{theorem}

\begin{theorem}\label{vertex}Let $p>3D$. Let $k=1,2,\cdots,d+e-1$ be such that $V_k$ contains exactly one pair. Then
$$\text{ord}_q(c_{bk})\geq p_{[-e,d]}(k)$$ with equality holding if
and only if $\bar{H}_k(\vec{a})\neq0$.
\end{theorem}

From now on, we suppose that $q=p^b$, and let $\zeta$ be a primitive
$(q-1)$-th roots of unity.

\begin{definition}We define
the matrix $G=(G_{(i,u),(j,w)})_{-e\leq i,j\leq d-1,0\leq u,w\leq
b-1}$ by
$$(\zeta^w)^{\sigma^{-1}}F_{ij}^{\sigma^{-1}}=\sum\limits_{u=0}^{b-1}G_{(i,u),(j,w)}\zeta^u.$$\end{definition}

\begin{lemma}We have
$$\Psi_p(\zeta^wx^j)\equiv\sum\limits_{i=-e}^{d-1}\sum\limits_{u=0}^{b-1}G_{(i,u),(j,w)}\zeta^ux^i (\mod \partial L).$$
That is, $G$ is the matrix of $\bar{\Psi}_p$ with respect to the
basis over $\mathbb{Q}_p(\pi^{1/D})$ represented by
$$\zeta^ux^i,\ -e\leq i\leq d-1, 0\leq u\leq b-1.$$\end{lemma}
{\it Proof. } Recall that
$$\psi_p\circ E_{\vec{a}}(x^j)\equiv\sum\limits_{i=-e}^{d-1}F_{ij}x^i(\mod \partial L).$$ So
$$\Psi_p(\zeta^wx^j)\equiv(\zeta^w)^{\sigma^{-1}}\sum\limits_{i=-e}^{d-1}F_{ij}^{\sigma^{-1}}x^i(\mod \partial L).$$
By definition,
$$(\zeta^w)^{\sigma^{-1}}F_{ij}^{\sigma^{-1}}=\sum\limits_{u=0}^{b-1}G_{(i,u),(j,w)}\zeta^u.$$
The lemma  now follows.

\begin{corollary}We have
$$\det(1-s\bar{\Psi}_p\mid L/(\partial L)\text{ \rm over
 }\mathbb{Q}_p(\pi^{1/D}))=\det(1-sG).$$
In
particular,
$$c_{bk}=\sum\limits_{T}\det((G_{(i,u),(j,w)})_{(i,u),(j,w)\in
T}),$$ where $T$ runs over subsets of
$$\{(i,u)\mid -e\leq i\leq d-1, 0\leq u\leq b-1\}$$ with cardinality
$bk$. \end{corollary}

\begin{lemma}\label{inequality}Let $T_1$ and $T_2$ be two finite sets with equal cardinality.
Let $g_1$ and $g_2$ be real-valued functions on $T_1$ and $T_2$
respectively. Suppose that $g_1$ and $g_2$ agree on $T_1\cap T_2$,
and that $g_2(t_2)\geq g_1(t_1)$ for $t_2\in T_2\setminus T_1$ and
$t_1\in T_1\setminus T_2$. Then $$\sum\limits_{t\in T_2}g_2(t)\geq
\sum\limits_{t\in T_1}g_1(t).$$ Moreover, if $g_2(t_2)>g_1(t_1)$ for
$t_2\in T_2\setminus T_1$ and $t_1\in T_1\setminus T_2$, then the
equality holds if and only if $T_1=T_2$.\end{lemma} {\it Proof.
}Obvious.

We are now ready to prove Theorem \ref{non-vertex}.

{\it Proof of Theorem \ref{non-vertex}. } It suffices to show that,
for any subset $T$ of
$$\{(i,u)\mid -e\leq i\leq d-1, 0\leq u\leq b-1\}$$ with cardinality
$bk$, and any permutation $\tau$ of $T$, we have
$$\text{ord}_{\pi}(\prod_{(i,u)\in T}G_{(i,u),\tau(i,u)})
\geq b(p-1)p_{[-e,d]}(k).$$

Let $V_k=\{(m-1,n+1),(m,n)\}$. Then $\frac{n+1}{d}=\frac{m}{e}$.
Moreover, the cardinality of the set $\{1\leq i\leq m-1\mid pi\equiv
m(\mod e)\}$ is equal to that of $\{1\leq i\leq n\mid pi\equiv
n+1(\mod d)\}$. Without loss of generality, we assume that both of
them are of cardinality $1$.
 Then
$$(p-1)p_{[-e,d]}(k)=\sum\limits_{i=1}^{n}\lceil\frac{pi-n}{d}\rceil
+\sum\limits_{i=1}^{m-1}\lceil\frac{pi-m+1}{e}\rceil+\lceil\frac{(p-1)m}{e}\rceil-1.$$

Note that
$$\text{ord}_{\pi}(G_{(i,u),\tau(i,u)})=\text{ord}_{\pi}(F_{i,\tau(i)}).$$
So, if $i>0$, then
$$\text{ord}_{\pi}(G_{(i,u),\tau(i,u)}) \geq\left\{
                                              \begin{array}{ll}
                                                \lceil\frac{pi-n}{d}\rceil, & \hbox{ } \tau(i)\leq n,\\
                                                \lceil\frac{pi-n}{d}\rceil-1, & \hbox{ } \tau(i)> n,\\
                                                \lceil\frac{p(n+1)}{d}\rceil+1, & \hbox{
}i>n+1.
                                              \end{array}
                                            \right.$$

Similarly, if $i<0$, then
$$\text{ord}_{\pi}(G_{(i,u),\tau(i,u)})
\geq\left\{
                                              \begin{array}{ll}
                                                \lceil\frac{-pi-m+1}{e}\rceil, & \hbox{ } \tau(i)\geq -m+1,\\
                                                \lceil\frac{-pi-m+1}{e}\rceil-1, & \hbox{ } \tau(i)\leq-m,\\
                                                \lceil\frac{pm}{e}\rceil+1, & \hbox{
}i<-m.
                                              \end{array}
                                            \right.$$

Therefore
$$\text{ord}_{\pi}(\prod_{(i,u)\in T}G_{(i,u),\tau(i,u)})
\geq\sum\limits_{(i,u)\in T:1\leq i\leq
n}\lceil\frac{pi-n}{d}\rceil+\sum\limits_{(-i,u)\in T:1\leq i\leq
m-1}\lceil\frac{pi-m+1}{e}\rceil$$$$+\lceil\frac{(p-1)m}{e}\rceil\sum\limits_{(i,u)\in
T:i=n+1\textrm{ or }i=-m}1+\sum\limits_{(i,u)\in T:i>n+1\textrm{ or
}i<-m}\lceil\frac{pm}{e}\rceil$$$$+\sum\limits_{(i,u)\in
T:i>n+1\textrm{ or }i<-m}1 -\sum\limits_{(i,u)\in
T:\tau(i)>n\textrm{ or }\leq-m}1$$
$$\geq\sum\limits_{(i,u)\in T:1\leq i\leq
n}\lceil\frac{pi-n}{d}\rceil+\sum\limits_{(-i,u)\in T:1\leq i\leq
m-1}\lceil\frac{pi-m+1}{e}\rceil$$$$+(\lceil\frac{(p-1)m}{e}\rceil-1)\sum\limits_{(i,u)\in
T:i=n+1\textrm{ or }i=-m}1+\sum\limits_{(i,u)\in T:i>n+1\textrm{ or
}i<-m}\lceil\frac{pm}{e}\rceil.$$ By Lemma \ref{inequality}, we have
$$\text{ord}_{\pi}(\prod_{(i,u)\in T}G_{(i,u),\tau(i,u)}) \geq
b(\sum\limits_{i=1}^{n}\lceil\frac{pi-n}{d}\rceil
+\sum\limits_{i=1}^{m-1}\lceil\frac{pi-m+1}{e}\rceil+\lceil\frac{(p-1)m}{e}\rceil-1).$$The
proof is completed.

It remains to prove Theorem \ref{vertex}.
\begin{lemma}Let $p>3D$. Let $k=1,2,\cdots,d+e-1$ be such that $V_k=\{(m,n)\}$. Then
$$c_{bk}=\det((G_{(i,u),(j,w)})_{-m\leq i,j\leq n,0\leq u,w\leq b-1})+O(\pi^{b(p-1)p_{[-e,d]}(k)+\frac{1}{D}}).$$
\end{lemma}

 {\it Proof. }It suffices to show that, for any
subset $T$ of
$$\{(i,u)\mid -e\leq i\leq d-1, 0\leq u\leq b-1\}$$ with cardinality
$bk$ which is different from
$\{-m,\cdots,n\}\times\{0,\cdots,b-1\}$, and any permutation $\tau$
of $T$, we have
$$\text{ord}_{\pi}(\prod_{(i,u)\in T}G_{(i,u),\tau((i,u))})>b(p-1)p_{[-e,d]}(k).$$

First we suppose that $I_k=\{(m,n)\}$. Note that, if $i>0$, then
$$\text{ord}_{\pi}(G_{(i,u),\tau(i,u)}) \geq\left\{
                                              \begin{array}{ll}
                                                \lceil\frac{pi-n}{d}\rceil, & \hbox{ } \tau(i)\leq n,\\
                                                \lceil\frac{pi-n}{d}\rceil-1, & \hbox{ } \tau(i)> n,\\
                                                \lceil\frac{pn}{d}\rceil+\frac{p}{d}-2, & \hbox{
}i>n.
                                              \end{array}
                                            \right.$$

Similarly, if $i<0$, then
$$\text{ord}_{\pi}(G_{(i,u),\tau(i,u)})
\geq\left\{
                                              \begin{array}{ll}
                                                \lceil\frac{-pi-m}{e}\rceil, & \hbox{ } \tau(i)\geq -m,\\
                                                \lceil\frac{-pi-m}{e}\rceil-1, & \hbox{ } \tau(i)<-m,\\
                                                \lceil\frac{pm}{e}\rceil+\frac{p}{e}-2, & \hbox{
}i<-m.
                                              \end{array}
                                            \right.$$
So
$$\text{ord}_{\pi}(\prod_{(i,u)\in T}G_{(i,u),\tau(i,u)})
\geq\sum\limits_{(i,u)\in T:1\leq i\leq
n}\lceil\frac{pi-n}{d}\rceil+\sum\limits_{(-i,u)\in T:1\leq i\leq
m}\lceil\frac{pi-m}{e}\rceil$$$$+\sum\limits_{(i,u)\in T:i>
n}\lceil\frac{pn}{d}\rceil+\sum\limits_{(-i,u)\in T: i>
m}\lceil\frac{pm}{e}\rceil$$$$+\sum\limits_{(i,u)\in T:i>n\textrm{
or }i<-m}(\frac{p}{D}-2) -\sum\limits_{(i,u)\in T:\tau(i)>n\textrm{
or }<-m}1$$$$>\sum\limits_{(i,u)\in T:1\leq i\leq
n}\lceil\frac{pi-n}{d}\rceil+\sum\limits_{(-i,u)\in T:1\leq i\leq
m}\lceil\frac{pi-m}{e}\rceil$$$$+\sum\limits_{(i,u)\in T:i>
n}\lceil\frac{pn}{d}\rceil+\sum\limits_{(-i,u)\in T: i>
m}\lceil\frac{pm}{e}\rceil$$ By Lemma \ref{inequality}, we have
$$\text{ord}_{\pi}(\prod_{(i,u)\in T}G_{(i,u),\tau(i,u)})>
b(p-1)p_{[-e,d]}(k).$$

Secondly, we suppose that $I_k$ contains two pairs. Without loss of
generality, we may assume that $I_k=\{(m,n),(m+1,n-1)\}$. Then
$\frac{m+1}{e}=\frac{n}{d}$, $$pi\not\equiv m+1(\mod e),\ 1\leq
i\leq m,$$ and there is exactly one $1\leq i\leq n-1$ such that
$$pi\equiv n(\mod d).$$
So
$$(p-1)p_{[-e,d]}(k)=\sum\limits_{i=1}^{n-1}\lceil\frac{pi-n+1}{d}\rceil
+\sum\limits_{i=1}^{m}\lceil\frac{pi-m-1}{e}\rceil+\lceil\frac{(p-1)n}{d}\rceil-1.$$
Note that, if $i>0$, then
$$\text{ord}_{\pi}(G_{(i,u),\tau(i,u)}) \geq\left\{
                                              \begin{array}{ll}
                                                \lceil\frac{pi-n+1}{d}\rceil, & \hbox{ } \tau(i)\leq n-1,\\
                                                \lceil\frac{pi-n+1}{d}\rceil-1, & \hbox{ } \tau(i)\geq n,\\
                                               \lceil\frac{pn}{d}\rceil+\frac{p}{d}-2, & \hbox{
}i>n.
                                              \end{array}
                                            \right.$$

Similarly, if $i<0$, then
$$\text{ord}_{\pi}(G_{(i,u),\tau(i,u)})
\geq\left\{
                                              \begin{array}{ll}
                                                \lceil\frac{-pi-m-1}{e}\rceil, & \hbox{ } \tau(i)\geq -m-1,\\
                                                \lceil\frac{-pi-m-1}{e}\rceil-1, & \hbox{ } \tau(i)<-m-1,\\
                                                \lceil\frac{pn}{d}\rceil+\frac{p}{e}-2, & \hbox{
}i<-m-1.
                                              \end{array}
                                            \right.$$

So
$$\text{ord}_{\pi}(\prod_{(i,u)\in T}G_{(i,u),\tau(i,u)})
\geq\sum\limits_{(i,u)\in T:1\leq i<
n}\lceil\frac{pi-n+1}{d}\rceil+\sum\limits_{(-i,u)\in T:1\leq i\leq
m}\lceil\frac{pi-m-1}{e}\rceil$$$$+\lceil\frac{(p-1)n}{d}\rceil\sum\limits_{(i,u)\in
T:i=n\textrm{ or }i=-m-1}1+\sum\limits_{(i,u)\in T:i>n\textrm{ or
}i<-m-1}\lceil\frac{pn}{d}\rceil$$$$+\sum\limits_{(i,u)\in
T:i>n\textrm{ or }i<-m-1}(\frac{p}{D}-2) -\sum\limits_{(i,u)\in
T:\tau(i)\geq n\textrm{ or }<-m-1}1.$$

If $\{(i,u)\in T:i>n\textrm{ or }i<-m-1\}\neq\emptyset$, then
$$\text{ord}_{\pi}(\prod_{(i,u)\in T}G_{(i,u),\tau(i,u)})>\sum\limits_{(i,u)\in T:1\leq i<
n}\lceil\frac{pi-n+1}{d}\rceil+\sum\limits_{(-i,u)\in T:1\leq i\leq
m}\lceil\frac{pi-m-1}{e}\rceil$$$$+(\lceil\frac{(p-1)n}{d}\rceil-1)\sum\limits_{(i,u)\in
T:i=n}1+\lceil\frac{(p-1)n}{d}\rceil\sum\limits_{(i,u)\in
T:i=-m-1}1+\sum\limits_{(i,u)\in T:i>n\textrm{ or
}i<-m-1}\lceil\frac{pn}{d}\rceil.$$

By Lemma \ref{inequality}, we have
$$\text{ord}_{\pi}(\prod_{(i,u)\in
T}G_{(i,u),\tau(i,u)})>\sum\limits_{i=1}^{n-1}\lceil\frac{pi-n+1}{d}\rceil
+\sum\limits_{i=1}^{m}\lceil\frac{pi-m-1}{e}\rceil+\lceil\frac{(p-1)n}{d}\rceil-1.$$

If $\{(i,u)\in T:i>n\textrm{ or }i<-m-1\}=\emptyset$, then
$$\text{ord}_{\pi}(\prod_{(i,u)\in T}G_{(i,u),\tau(i,u)})\geq\sum\limits_{(i,u)\in T:1\leq i<
n}\lceil\frac{pi-n+1}{d}\rceil+\sum\limits_{(-i,u)\in T:1\leq i\leq
m}\lceil\frac{pi-m-1}{e}\rceil$$$$+(\lceil\frac{(p-1)n}{d}\rceil-1)\sum\limits_{(i,u)\in
T:i=n}1+\lceil\frac{(p-1)n}{d}\rceil\sum\limits_{(i,u)\in
T:i=-m-1}1.$$

By Lemma \ref{inequality}, we also have
$$\text{ord}_{\pi}(\prod_{(i,u)\in
T}G_{(i,u),\tau(i,u)})>\sum\limits_{i=1}^{n-1}\lceil\frac{pi-n+1}{d}\rceil
+\sum\limits_{i=1}^{m}\lceil\frac{pi-m-1}{e}\rceil+\lceil\frac{(p-1)n}{d}\rceil-1.$$

The proof is completed.

\begin{definition}We write $\alpha\sim\beta$ to mean that $\alpha=u\beta$ for some $p$-adic unit $u$.
\end{definition}

\begin{theorem}\label{semi-final}Let $p>3D$. Let $k=1,2,\cdots,d+e-1$ be such that $V_k=\{(m,n)\}$.
Then
$$c_{bk}\sim\det((F_{ij})_{-m\leq i,j\leq n})^b+O(\pi^{b(p-1)p_{[-e,d]}(k)+\frac{1}{D}}).$$
\end{theorem}
{\it Proof. }It suffices to show that
$$\det((G_{(i,u),(j,w)})_{-m\leq i,j\leq n,0\leq u,w\leq b-1})\sim\det((F_{ij})_{-m\leq i,j\leq n})^b.$$
Let $V=\oplus_{i=-m}^n\mathbb{Q}_q(\pi^{1/D})e_i$ be a
$k$-dimensional vector space over $\mathbb{Q}_q(\pi^{1/D})$with
standard basis $e_{-m},\cdots,e_n$. Let $F=(F_{ij})_{-m\leq i,j\leq
n}$ act on it in the standard way, and let $\sigma$ act on it
coordinate-wise. Then
$$\sigma^{-1}\circ F(\zeta^we_j)=(\zeta^w)^{\sigma^{-1}}\sum\limits_{i=-m}^nF_{ij}^{\sigma^{-1}}e_i.$$
Therefore, $G$ is the matrix of $\sigma^{-1}\circ F$ on $V$ with
respect to the basis over $\mathbb{Q}_p(\pi^{1/D})$:
$$\zeta^ue_i,\ -m\leq i\leq n, 0\leq u\leq b-1.$$
As $\sigma$ is just a re-ordering of the basis, we have
$$\det((G_{(i,u),(j,w)})_{-m\leq i,j\leq n,0\leq u,w\leq b-1})\sim\det((F_{ij})_{-m\leq i,j\leq n})^b.$$
The theorem is proved.

\begin{lemma}\label{vetex-main-term}Let $p>3D$. Let $k=1,2,\cdots,d+e-1$ be such that $V_k=\{(m,n)\}$.
Then
$$\det((F_{ij})_{-m\leq i,j\leq n})=\sum\limits_{\tau\in S_k}\text{sgn}(\tau)\prod_{i=-m}^nF_{i,\tau(i)}
+O(\pi^{(p-1)p_{[-e,d]}(k)+1/D}).$$
\end{lemma}
{\it Proof. } For $j\leq n$, we have
$$\lceil\frac{pi-j}{d}\rceil=\lceil\frac{pi-n+(n-j)}{d}\rceil\geq\lceil\frac{pi-n}{d}\rceil$$
with equality holding if and only if
$$j\geq n-d\{-\frac{pi-n}{d}\}.
$$
Similarly, for $j\geq-m$, we have
$$\lceil\frac{-pi+j}{e}\rceil=\lceil\frac{-pi-m+(m+j)}{e}\rceil
\geq\lceil\frac{-pi-m}{e}\rceil
$$
with equality holding if and only if
$$j\leq -m+e\{\frac{pi+m}{e}\}.
$$
So, if $\tau\not\in S_k$ is a permutation of
$\{-m,-(m-1),\cdots,n\}$, then
$$\text{ord}_{\pi}(\prod_{i=-m}^nF_{i,\tau(i)})\geq
\sum\limits_{i=-m}^n\lceil\deg(pi-\tau(i))\rceil
$$$$\geq1+\sum\limits_{i=1}^n\lceil\frac{pi-n}{d}\rceil+
\sum\limits_{i=1}^m\lceil\frac{pi-m}{e}\rceil.$$
Hence$$\det((F_{ij})_{-m\leq i,j\leq n})=\sum\limits_{\tau\in
S_k}\text{sgn}(\tau)\prod_{i=-m}^nF_{i,\tau(i)}+O(\pi^{(p-1)p_{[-e,d]}(k)+1/D}).$$
The lemma is proved.

We are now ready to prove Theorem \ref{vertex}. By the above lemmas,
it suffices to prove the following theorem.
\begin{theorem}\label{final-approximation}Let $p>3D$. Let $k=1,2,\cdots,d+e-1$ be such that $V_k=\{(m,n)\}$.
Then
$$\det((F_{ij})_{-m\leq i,j\leq n})=\pi^{(p-1)p_{[-e,d]}(k)}\hat{a}_d^{u_k}
\hat{a}_{-e}^{v_k}H_k(\vec{\hat{a}})+O(\pi^{(p-1)p_{[-e,d]}(k)+1/D}),$$
where $\vec{\hat{a}}=(\hat{a}_{-e},\cdots,\hat{a}_d)$, and $u_k,v_k$
are integers depending on $k$.
\end{theorem}
{\it Proof. } By Lemmas \ref{vetex-main-term} and
\ref{frobenius-main-term}, we have
$$\det((F_{ij})_{-m\leq i,j\leq n})=\sum\limits_{\tau\in
S_k}\text{sgn}(\tau)\prod_{i=-m}^n\gamma_{pi-\tau(i)}
+O(\pi^{(p-1)p_{[-e,d]}(k)+1/D}).$$ By Lemmas
\ref{artin-hasse-plus-main-term} and
\ref{artin-hasse-minus-main-term}, we have
$$\gamma_{pi-\tau(i)}=\left\{
                        \begin{array}{ll}
                          \pi^{\lceil\frac{pi-n}{d}\rceil}\lambda_{\lfloor\frac{pi-\tau(i)}{d}\rfloor}
\lambda_{\lceil\{\frac{pi-\tau(i)}{d}\}\rceil}
\hat{a}_d^{\lceil\frac{pi-n}{d}\rceil-1} \hat{a}_{r_i-\tau(i)}
+O(\pi^{\lceil\frac{pi-n}{d}\rceil}), & \hbox{ } i>0\\
                          \pi^{\lceil\frac{-pi-m}{e}\rceil}\lambda_{\lfloor\frac{-pi+\tau(i)}{e}\rfloor}
\lambda_{\lceil\{\frac{-pi+\tau(i)}{e}\}\rceil}\hat{a}_{-e}^{\lceil\frac{-pi-m}{e}\rceil-1}
\hat{a}_{-r_i-\tau(i)} +O(\pi^{\lceil\frac{-pi-m}{e}\rceil+1}), &
\hbox{ }i<0.
                        \end{array}
                      \right.$$
The theorem now follows.


\end{document}